\documentclass[a4paper, 11pt, reqno, final]{amsart}

\usepackage{amsmath,amssymb,amsthm}
\usepackage{graphicx}
\usepackage[matrix,arrow]{xy}

\numberwithin{equation}{section}
\allowdisplaybreaks

\theoremstyle{definition}
 \newtheorem{thm}{Theorem}[section]
 \newtheorem{prp}[thm]{Proposition}
 \newtheorem{lem}[thm]{Lemma}
 
 \newtheorem{dfn}[thm]{Definition}
 \newtheorem{fct}[thm]{Fact}
 \newtheorem{rmk}[thm]{Remark}
 
 \newtheorem*{ack}{Acknowledgements}
 \newtheorem*{not*}{Notation}
 \newtheorem*{rmk*}{Remark}
 \newtheorem*{thm*}{Theorem}

\newcommand{\bbC}{\mathbb{C}}

\newcommand{\bbQ}{\mathbb{Q}}
\newcommand{\bbZ}{\mathbb{Z}}
\newcommand{\calF}{\mathcal{F}}

\newcommand{\calH}{\mathcal{H}}
\newcommand{\calL}{\mathcal{L}}
\newcommand{\calP}{\mathcal{P}}
\newcommand{\calW}{\mathcal{W}}
\newcommand{\frkg}{\mathfrak{g}}
\newcommand{\frkn}{\mathfrak{n}}
\newcommand{\frksl}{\mathfrak{sl}}

\newcommand{\ep}{\epsilon}
\newcommand{\Vir}{\mathrm{Vir}}
\newcommand{\SU}{\mathrm{SU}}
\newcommand{\rk}{\mathrm{rank}}

\newcommand{\seteq}{\mathbin{:=}}
\newcommand{\ket}[1]{\left|  #1  \right>}
\newcommand{\bra}[1]{\left< #1 \right|}
\newcommand{\prd}[1]{\left<\right.\hskip-0.25em #1 \hskip-0.25em \left.\right>}
\newcommand{\nop}{\genfrac{}{}{0pt}{1}{\circ}{\circ}}
\newcommand{\no}[1]{\nop #1 \nop}
\newcommand{\vect}[1]{\overrightarrow{#1}}

\makeatletter
\def\Res{\mathop{\operator@font Res}}
\def\Ind{\mathop{\operator@font Ind}\nolimits}
\def\End{\mathop{\operator@font End}}
\makeatother

\title[Whittaker vector of Virasoro algebra via Jack polynomial]{Whittaker vectors of the Virasoro algebra in terms of Jack symmetric polynomial}
\date{March 4, 2010; revised March 8, 2011}
\author{Shintarou Yanagida}

\address{Kobe University, Department of Mathematics, Rokko, Kobe 657-8501, Japan}
\email{yanagida@math.kobe-u.ac.jp}

\begin{document}

\begin{abstract}
We give an explicit formula of Whittaker vector for Virasoro algebra 
in terms of the Jack symmetric functions.
Our fundamental tools are the Feigin-Fuchs bosonization and 
the split expression of the Calogero-Sutherland model given by 
Awata-Matsuo-Odake-Shiraishi.
\end{abstract}

\maketitle

\section{Introduction}\label{sect:intro}

In \cite{AGT:2010} a remarkable proposal, now called the AGT conjecture, 
was given on the relation between the Liouville theory conformal blocks 
and the Nekrasov partition function.
Among the related investigations, 
Gaiotto proposed several degenerated versions of the AGT conjecture 
in \cite{G:2009}.
In that paper, he conjectured that the inner product of a certain element 
in the Verma module of Virasoro algebra coincides 
with the Nekrasov partition function 
for the four dimensional $\mathcal{N}=2$ pure gauge theory \cite{N:2003}.
Actually, the element considered is a kind of Whittaker vector 
in the Verma module of the Virasoro algebra.

Whittaker vectors and Whittaker modules are important gadgets 
in the representation theory since its emergence in the study of 
finite dimensional Lie algebras \cite{K:1978}.
Although numerous analogues and generalisations have been proposed for 
other algebras, such as affine algebras and quantum groups,
not so many investigations have been given for
the Whittaker vectors of the Virasoro algebra.
A general theory on the properties of Whittaker modules 
for the Virasoro algebra was recently given in \cite{OW:2009}.

In this paper we give an explicit expression of the Whittaker vector for 
the Verma module of Virasoro algebra in terms of Jack symmetric functions 
\cite[VI \S 10]{M:1995:book}.
We use the Feigin-Fuchs bosonization \cite{FF:1982} 
to identify the Verma module and the ring of symmetric function, 
and then utilise the split expression of the Calogero-Sutherland Hamiltonian \cite{SKAO:1996} to derive an recursion relation on the coefficients 
of the Whittaker vector in its expansion 
with respect to Jack symmetric functions.

Our result is related to a conjecture given 
by Awata and Yamada in \cite{AY:2009}.
They proposed the five-dimensional AGT conjecture 
 for pure $\SU(2)$ gauge theory using  the deformed Virasoro algebra, 
and as a related topic, they also proposed a conjectural formula 
on the explicit form of the deformed Gaiotto state 
in terms of Macdonald symmetric functions \cite[(3.18)]{AY:2009}.
Our formula is the non-deformed Virasoro, or four-dimensional, counterpart 
of their conjectural formula.

The motivation of our study also comes from the work \cite{MY:1995},
where singular vectors of the Virasoro algebra are 
expressed by Jack polynomials. 

Before presenting the detail of the main statement,
we need to prepare several notations on Virasoro algebra,
symmetric functions and some combinatorics.
The main theorem will be given in \S \ref{subsec:mainthm}.

\subsection{Partitions}\label{subsec:partition}

Throughout in this paper, notations of partitions follow \cite{M:1995:book}. 
For the positive integer $n$, a partition $\lambda$ of $n$ is a (finite) 
sequence of positive integers $\lambda=(\lambda_1,\lambda_2,\ldots)$ such that 
$\lambda_1\ge\lambda_2\ge\cdots$ and $\lambda_1+\lambda_2+\cdots=n$.
The symbol $\lambda \vdash n$ means that $\lambda$ is a partition of $n$. 
For a general partition we also define $|\lambda|\seteq\sum_{i} \lambda_i$.
The number $\ell(\lambda)$ is defined to be 
the length of the sequence $\lambda$.
The conjugate partition of $\lambda$ is denoted by $\lambda'$.

We also consider the empty sequence $\emptyset$ 
as the unique partition of the number $0$.

In addition we denote by $\calP$ 
the set of all the partitions of natural numbers 
including the empty partition $\emptyset$.
So that we have
\begin{align*}
\calP=\{\emptyset,(1),(2),(1,1),(3),(2,1),(1,1,1),\ldots\}.
\end{align*}
As usual, 
$p(n)\seteq\#\{\lambda\in\calP\mid |\lambda|=n\}
=\#\{\lambda\mid \lambda\vdash n\}$ 
denotes the number of partitions of $n$.

In the main text, we sometimes use the dominance semi-ordering 
on the partitions: 
$\lambda\ge\mu$ if and only if
$|\lambda|=|\mu|$ and 
$\sum_{k=1}^i \lambda_k \ge \sum_{k=1}^i \mu_k$ ($i=1,2,\ldots$).

\begin{figure}[htbp]
\unitlength 0.1in
\begin{center}
\begin{picture}(  8.0000,  12.0000)(  4.0000,-15.0000)
\special{pn 8}%
\special{pa 400 400}%
\special{pa 600 400}%
\special{pa 600 600}%
\special{pa 400 600}%
\special{pa 400 400}%
\special{fp}%
\special{pn 8}%
\special{pa 600 400}%
\special{pa 800 400}%
\special{pa 800 600}%
\special{pa 600 600}%
\special{pa 600 400}%
\special{fp}%
\special{pn 8}%
\special{pa 800 400}%
\special{pa 1000 400}%
\special{pa 1000 600}%
\special{pa 800 600}%
\special{pa 800 400}%
\special{fp}%
\special{pn 8}%
\special{pa 1000 400}%
\special{pa 1200 400}%
\special{pa 1200 600}%
\special{pa 1000 600}%
\special{pa 1000 400}%
\special{fp}%
\special{pn 8}%
\special{pa 400 600}%
\special{pa 600 600}%
\special{pa 600 800}%
\special{pa 400 800}%
\special{pa 400 600}%
\special{fp}%
\special{pn 8}%
\special{pa 600 600}%
\special{pa 800 600}%
\special{pa 800 800}%
\special{pa 600 800}%
\special{pa 600 600}%
\special{fp}%
\special{pn 8}%
\special{pa 800 600}%
\special{pa 1000 600}%
\special{pa 1000 800}%
\special{pa 800 800}%
\special{pa 800 600}%
\special{fp}%
\special{pn 8}%
\special{pa 1000 600}%
\special{pa 1200 600}%
\special{pa 1200 800}%
\special{pa 1000 800}%
\special{pa 1000 600}%
\special{fp}%
\special{pn 8}%
\special{pa 400 800}%
\special{pa 600 800}%
\special{pa 600 1000}%
\special{pa 400 1000}%
\special{pa 400 800}%
\special{fp}%
\special{pn 8}%
\special{pa 600 800}%
\special{pa 800 800}%
\special{pa 800 1000}%
\special{pa 600 1000}%
\special{pa 600 800}%
\special{fp}%
\special{pn 8}%
\special{pa 400 1000}%
\special{pa 600 1000}%
\special{pa 600 1200}%
\special{pa 400 1200}%
\special{pa 400 1000}%
\special{fp}%
\special{pn 8}%
\special{pa 400 1200}%
\special{pa 600 1200}%
\special{pa 600 1400}%
\special{pa 400 1400}%
\special{pa 400 1200}%
\special{fp}%
\special{pn 8}%
\special{pa 400 1400}%
\special{pa 600 1400}%
\special{pa 600 1600}%
\special{pa 400 1600}%
\special{pa 400 1400}%
\special{fp}%
\put(3.0, -3.0){\vector(1,0){4}}
\put(3.0, -3.0){\vector(0,-1){4}}
\put(7.6, -3.4){\makebox(0,0)[lb]{$j$}}%
\put(2.8, -8.5){\makebox(0,0)[lb]{$i$}}%
\end{picture}%
\end{center}
\caption{The Young diagram for $(4,4,2,1,1,1)$}
\label{fig:442111}
\end{figure}
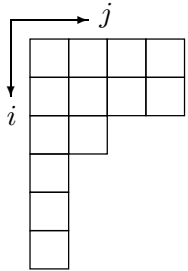

We also follow \cite{M:1995:book} for the convention of the Young diagram.
Moreover we will use the coordinate $(i,j)$ on the Young diagram 
defined as follows:
the first coordinate $i$ (the row index) increases as one goes downwards,
and the second coordinate $j$ (the column index) increases as one goes 
rightwards. 
For example, in Figure \ref{fig:442111} 
the left-top box has the coordinate $(1,1)$ and 
the left-bottom box has the coordinate $(6,1)$.
We will often identify a partition and its associated Young diagram.

Let us also use the notation $(i,j)\in\lambda$, 
which means that $i,j\in\bbZ_{\ge 1}$,
$1\le i\le \ell(\lambda)$ and $1\le j\le \lambda_i$.
On the Young diagram of $\lambda$ the symbol $(i,j)\in\lambda$ corresponds 
to the box located at the coordinate $(i,j)$.
In Figure \ref{fig:442111}, we have $(2,3)\in\lambda\seteq(4,4,2,1,1)$ 
but $(4,3)\notin\lambda$. 

\subsection{Virasoro algebra}\label{subsec:vir}

Let us fix notations on Virasoro algebra and its Verma module.
Let $c\in\bbC$ be a fixed complex number.
The Virasoro algebra $\Vir_c$ is a Lie algebra over $\bbC$ 
with central extension,
generated by $L_n$ ($n\in\bbZ$) with the relation
\begin{align}\label{eq:Vir}
[L_m,L_n]=(m-n)L_{m+n}+\dfrac{c}{12}m(m^2-1)\delta_{m+n,0}.
\end{align}
$\Vir_c$ has the triangular decomposition 
$\Vir_c=\Vir_{c,+}\oplus\Vir_{c,0}\oplus\Vir_{c,-}$ with 
$\Vir_{c,\pm} \seteq \oplus_{\pm n>0}\bbC L_n$ and
$\Vir_{c,0} \seteq \bbC L_0\oplus \bbC$.

Let $h$ be a complex number.
Let $\bbC_{h}$ be the one-dimensional representation 
of the subalgebra $\Vir_{c,\ge0} \seteq \Vir_{c,0}\oplus \Vir_{c,+}$,
where $\Vir_{c,+}$ acts trivially and $L_0$ acts as the multiplication by $h$.
Then one has the Verma module $M_h$ by
\begin{align*}
M_h\seteq\Ind_{\Vir_{c,\ge0}}^{\Vir_c}\bbC_{h}.
\end{align*}
Obeying the notation in physics literature, 
we denote by $\ket{h}$ a fixed basis of $\bbC_{h}$.
Then one has $\bbC_{h}=\bbC \ket{h}$ and $M_h=U(\Vir_c)\ket{h}$.

$M_h$ has the $L_0$-weight space decomposition: 
\begin{align}\label{eq:L0weightdecomp}
M_h=\bigoplus_{n\in\bbZ_{\ge 0}} M_{h,n},\quad \text{with}\quad  
M_{h,n}\seteq\{v\in M_h\mid L_0 v=(h+n)v \}.
\end{align}

A basis of $M_{h,n}$ can be described simply by partitions.
For a partition 
$\lambda=(\lambda_1,\lambda_2,\ldots,\lambda_k)$ of $n$ 
we define the abbreviation
\begin{align}\label{eq:L-lambda}
L_{-\lambda}\seteq L_{-\lambda_k} L_{-\lambda_{k-1}}\cdots L_{-\lambda_1}.
\end{align}
of the element of $U(\Vir_{c,-})$, the enveloping algebra 
of the subalgebra $\Vir_{c,-}$.
Then the set
\begin{align*}
\{L_{-\lambda} \ket{h} \mid \lambda \vdash n\},
\end{align*}
is a basis of $M_{h,n}$.

\subsection{Bosonization}

Next we recall the bosonization of the Virasoro algebra \cite{FF:1982}. 
Consider the Heisenberg algebra $\calH$ generated by $a_n$ ($n\in\bbZ$) 
with the relation
\begin{align*}
[a_m,a_n]=m\delta_{m+n,0}.
\end{align*}
Consider the correspondence
\begin{align}
L_n \mapsto \calL_n
 \seteq\dfrac{1}{2}\sum_{m\in\bbZ}\no{a_m a_{n-m}}-(n+1)\rho a_n,
\label{eq:FF}
\end{align}
where the symbol $\no{\ }$ means the normal ordering.
This correspondence determines a well-defined morphism 
\begin{align}\label{eq:FF:phi}
\varphi: U(\Vir_c) \to\widehat{U}(\calH).
\end{align}
Here $\widehat{U}(\calH)$ is 
the completion of the universal enveloping algebra $U(\calH)$ 
in the following sense \cite{FF:1996}. 
For $n\in\bbZ_{\ge0}$, 
let $I_n$ be the left ideal of the enveloping algebra $U(\calH)$ 
generated by all polynomials in $a_m$ ($m\in\bbZ_{\ge1})$ 
of degrees greater than or equal to $n$ 
(where we defined the degree by $\deg a_m\seteq m$).
Then we define 
\begin{align*}
\widehat{U}(\calH)\seteq\varprojlim_n \widehat{U}(\calH)/I_n.
\end{align*}

Next we recall the functorial correspondence of the representations. 
First let us define the Fock representation $\calF_\alpha$ of $\calH$.
$\calH$ has the triangular decomposition 
$\calH=\calH_{+}\oplus\calH_{0}\oplus\calH_{-}$
with $\calH_{\pm}\seteq \oplus_{\pm n\in\bbZ_{+}}\bbC a_n$ and 
$\calH_{0}\seteq \bbC a_0$.
Let $\bbC_\alpha=\bbC\ket{\alpha}_\calF$ 
be the one-dimensional representation of 
$\calH_{0}\oplus\calH_{+}$ with the action 
$a_0\ket{\alpha}_\calF=\alpha \ket{\alpha}_\calF$ and 
$a_n\ket{\alpha}_\calF=0$ ($n\in\bbZ_{\ge1}$).
Then the Fock space $\calF_\alpha$ is defined to be 
\begin{align*}
\calF_\alpha\seteq
\Ind_{\calH_{0}\oplus\calH_{-}}^{\calH}\bbC_\alpha
\end{align*}
It has the $a_0$-weight decomposition  
\begin{align}\label{eq:F:deg}
\calF_\alpha=\oplus_{n\ge 0}\calF_{\alpha,n},\quad
\calF_{\alpha,n}
\seteq\{w\in\calF_{\alpha} \mid a_0 w= (n+\alpha)w\}.
\end{align}
Each weight space $\calF_{\alpha,n}$ has a basis 
\begin{align}\label{eq:a:base}
\{a_{-\lambda}\ket{\alpha}_\calF \mid \lambda\vdash n \}
\end{align} 
with $a_{-\lambda}\seteq a_{-\lambda_k}\cdots a_{-\lambda_1}$ 
for a partition $\lambda=(\lambda_1,\ldots,\lambda_k)$.
Note also that 
the action of $\widehat{U}(\calH)$ on $\calF_\alpha$ is well-defined.

Similarly the dual Fock space $\calF^*_\alpha$ is defined to be 
$\Ind_{\calH_{0}\oplus\calH_{-}}^{\calH }\bbC_\alpha^*$,
where $\bbC_\alpha^*=\bbC\cdot{}_\calF\bra{\alpha}$ 
is the one-dimensional right representation of 
$\calH_{0}\oplus\calH_{-}$ with the action 
${}_\calF\bra{\alpha}a_0=\alpha \cdot {}_\calF\bra{\alpha}$ and 
${}_\calF\bra{\alpha} a_{-n}=0$ ($n\in\bbZ_{\ge1}$).

Then one has the bilinear form 
\begin{align*}
\cdot: \calF^*_\alpha \times \calF_\alpha \to \bbC
\end{align*}
defined by
\begin{align*}
&{}_\calF\bra{\alpha}\cdot \ket{\alpha}_\calF =1,\quad 
0 \cdot \ket{\alpha}_\calF={}_\calF\bra{\alpha}\cdot 0=0,\\
&{}_\calF\bra{\alpha}u u'\cdot\ket{\alpha}_\calF
={}_\calF\bra{\alpha}u \cdot u'\ket{\alpha}_\calF
={}_\calF\bra{\alpha}\cdot u u'\ket{\alpha}_\calF\ 
(u,u'\in \calH).
\end{align*}
As in the physics literature, we often omit the symbol $\cdot$ and 
simply write ${}_\calF\left<\alpha|\alpha\right>_\calF$, 
${}_\calF\bra{\alpha}u\ket{\alpha}_\calF$ and so on.

Now we can state the bosonization of representation: 
\eqref{eq:FF} is compatible with the map 
\begin{align}\label{eq:FF:psi}
\psi: M_h \to \calF_\alpha,\quad 
L_{-\lambda}\ket{h}\mapsto \calL_{-\lambda}\ket{\alpha}_\calF
\end{align}
with $\calL_{-\lambda}\seteq\calL_{-\lambda_1}\cdots\calL_{-\lambda_k}$ 
for $\lambda=(\lambda_1,\ldots,\lambda_k)\in\calP$ and 
\begin{align}\label{eq:FF:ch}
c=1-12\rho^2,\quad h=\dfrac{1}{2}\alpha(\alpha-2\rho).
\end{align}
In other words, we have 
\begin{align*}
\psi(x v)=\varphi(x)\psi(v)\quad (x\in \Vir_c,\ v\in M_{h})
\end{align*} 
under the parametrisation \eqref{eq:FF:ch} 
of highest weights.

\subsection{Fock space and symmetric functions}
\label{subsec:Fs}

The Fock space $\calF_\alpha$ is naturally identified 
with the space of symmetric functions. 
In this paper the term ``symmetric function" means the 
infinite-variable symmetric ``polynomial". 
To treat such an object rigorously, 
we follow the argument of \cite[\S I.2]{M:1995:book}.

Let us denote by $\Lambda_N$ the ring of $N$-variable symmetric polynomials 
over $\bbZ$, 
and by $\Lambda_N^d$ the space of homogeneous symmetric polynomials 
of degree $d$. 
The ring of symmetric functions $\Lambda$ is defined as the inverse limit 
of the $\Lambda_N$ in the category of graded rings 
(with respect to the gradation defined by the degree $d$). 
We denote by $\Lambda_K=\Lambda\otimes_\bbZ K$ 
the coefficient extension to a ring $K$.
Among several bases of $\Lambda$, the family of 
the power sum symmetric functions 
\begin{align*}
p_n=p_n(x)\seteq \sum_{i\in\bbZ_{\ge1}} x_i^n,\quad
p_\lambda\seteq p_{\lambda_1}\cdots p_{\lambda_k},
\end{align*}
plays an important role. 
It is known that $\{p_\lambda\mid \lambda \vdash d\}$ is a basis of 
$\Lambda_\bbQ^d$, the subspace of homogeneous symmetric functions 
of degree $d$. 

Now following \cite{AMOS:1995}, 
we define the isomorphism between the Fock space 
and the space of symmetric functions.
Let $\beta$ be a non-zero complex number and consider the next map between 
$\calF_\alpha$ and $\Lambda_{\bbC}$: 
\begin{align}
\begin{array}{c c c c}
\iota_\beta: & \calF_\alpha & \to & \Lambda_{\bbC},
\\ 
& \rotatebox{90}{$\in$} & & \rotatebox{90}{$\in$}
\\
&v&\mapsto 
&{\displaystyle
{}_\calF\bra{\alpha}
\exp\Big(\sqrt{\dfrac{\beta}{2}}\sum_{n=1}^\infty \dfrac{1}{n}p_n a_n\Big)v.}
\end{array}
\label{eq:iota}
\end{align}
Under this morphism, 
an element $a_{-\lambda}\ket{\alpha}_\calF$ of the base \eqref{eq:a:base} 
is mapped to 
\begin{align*}
\iota_\beta(a_{-\lambda}\ket{\alpha}_\calF)
=(\sqrt{\beta/2})^{\ell(\lambda)} p_\lambda(x).
\end{align*}
Since $\{a_{-\lambda}\ket{\alpha}_\calF\}$ is a basis of $\calF_\alpha$ 
and $\{p_\lambda\}$ is a basis of $\Lambda_\bbQ$, 
$\iota_\beta$ is an isomorphism.

\subsection{Jack symmetric function}
\label{subsec:jack}

Now we recall the definition of Jack symmetric function 
\cite[\S VI.10]{M:1995:book}. 
Let $b$ be an indeterminate 
\footnote{Our parameter $b$ is usually denoted by $\alpha$ in the literature, for example, in \cite{M:1995:book}. We avoid using $\alpha$ since it is already defined to be the highest weight of the Heisenberg Fock space $\calF_\alpha$.}
and define an inner product on 
$\Lambda_{\bbQ(b)}$ by 
\begin{align}\label{eq:inner}
\left<p_\lambda,p_\mu\right>_b \seteq
\delta_{\lambda,\mu}z_\lambda b^{\ell(\lambda)}.
\end{align}
Here the function $z_\lambda$ is given by:
\begin{align*}
z_{\lambda}\seteq\prod_{i\in\bbZ_{\ge1}}i^{m_i(\lambda)} m_i(\lambda) ! 
\quad \text{ with } \quad
m_i(\lambda)\seteq \#\{1\le i\le \ell(\lambda) \mid \lambda_j=i \}.
\end{align*}
Then the (monic) Jack symmetric function $P_{\lambda}^{(b)}$ is determined 
uniquely by the following two conditions:
\begin{description}
\item[(i)] 
It has an expansion via monomial symmetric function $m_\nu$ in the form
\begin{align*}
P_{\lambda}^{(b)}
=m_\lambda+\sum_{\mu<\lambda}c_{\lambda,\mu}(b)m_\mu.
\end{align*}
Here $c_{\lambda,\mu}(b)\in\bbQ(b)$ and 
the ordering $<$ among the partitions is the dominance semi-ordering.

\item[(ii)]
The family of Jack symmetric functions is an orthogonal basis 
of $\Lambda_{\bbQ(b)}$ with respect to $\left<\cdot,\cdot\right>_b$:
\begin{align*}
\prd{P_\lambda^{(b)},P_\mu^{(b)}}_b=0 \quad 
\text{ if } \lambda\neq\mu.
\end{align*}
\end{description}


\subsection{Main Theorem}\label{subsec:mainthm}

Finally we can state our main statement.

Consider the Verma module $M_{h}$ of the Virasoro algebra $\Vir_c$ 
with generic complex numbers $c$ and $h$.
Let $a$ be an arbitrary complex number, 
and let $v_G$ be an element of the Verma module $M_{h}$ such that 
\[L_1 v_G=a v_G,\quad L_n v_G=0\ (n\ge 2).\] 
Then $v_G$ exists uniquely up to scalar multiplication
(see Fact \ref{fct:G:ue}).

Introduce the complex numbers $\rho$, $\alpha$ and $\beta$ by the relations
\begin{align*}
c=1-12\rho^2,\quad
h=\dfrac{1}{12}\alpha(\alpha-2\rho),\quad
\rho=\dfrac{\beta^{1/2}-\beta^{-1/2}}{\sqrt{2}}.
\end{align*}
Then by the Feigin-Fuchs bosonization $\psi:M_{h}\to \calF_\alpha$ 
\eqref{eq:FF:psi} 
and the isomorphism 
$\iota_\beta:\calF_\alpha \to \Lambda_{\bbC}$ 
\eqref{eq:iota},
one has an element $\iota_\beta \circ \psi (v_{G})\in\Lambda_\bbC$.

\begin{thm*}
We have
\begin{align}\label{eq:thm:expand}
\iota_\beta \circ \psi (v_{G})=
\sum_{\lambda\in\calP} 
 a^{|\lambda|} c_\lambda(\alpha,\beta) P_{\lambda}^{(\beta^{-1})},
\end{align}
where $\lambda$ runs over all the partitions and 
the coefficient $c_\lambda(\alpha,\beta)$ is given by 
\begin{align}\label{eq:thm}
\begin{split}
&c_\lambda(\alpha,\beta)\\
&=\prod_{(i,j)\in\lambda}\dfrac{1}{\lambda_i-j+1+\beta(\lambda_j'-i)}
 \prod_{\substack{(i,j)\in\lambda \\ (i,j)\neq (1,1) }}
 \dfrac{\beta}{(j+1)+\sqrt{2}\beta^{1/2} \alpha -(i+1)\beta}.
\end{split}
\end{align}
(See \S \ref{subsec:partition} for the symbol ``$(i,j)\in\lambda$".)
\end{thm*}

The proof of this theorem will be given in \S \ref{subsec:coh}.

In the main theorem above, the element $v_G$ is the Whittaker vector 
associated to the degenerate Lie algebra homomorphism 
$\eta:\Vir_{c,+}\to \bbC$, that is, $\eta(L_2)=0$.  
We shall call this element by ``Gaiotto state", following \cite{AY:2009}.
A general theory of Whittaker vectors usually assumes the non-degeneracy 
of the homomorphism $\eta$, i.e., $\eta(L_1)\neq0$ and $\eta(L_2)\neq0$.
This non-degenerate case will be treated in Proposition \ref{prp:rec2},
although there seem no factored expressions for the coefficients 
as \eqref{eq:thm}.

The content of this paper is as follows.
In \S \ref{sec:vir} we recall the split expression 
of the Calogero-Sutherland Hamiltonian,
which is a key point in our proof.
In \S \ref{sec:whit} we investigate the Whittaker vectors 
in terms of symmetric functions.
The main theorem will be proved \S \ref{subsec:coh}, 
using some combinatorial identities shown in \S \ref{sec:id}.
The Whittaker vector with respect to the non-degenerate homomorphism 
will be treated in \S \ref{subsec:wit}.
In the final \S \ref{sec:rmk} we give some remarks on possible 
generalisations and the related works.
We also added Appendix \ref{sec:AGT} concerning the AGT relation 
and its connection to our argument.

\section{Preliminaries on Jack symmetric functions and bosonized Calogero-Sutherland Hamiltonian}
\label{sec:vir}

This section is a preliminary for the proof of the main theorem.
We need the following Definition \ref{dfn:fE} and Proposition \ref{prp:split}:

\begin{dfn}\label{dfn:fE}
(1)
Let $\lambda$ be a partition and $b,\beta$ be generic complex numbers.
Define $f_{\lambda}^{(b,\beta)} \in \calF_{\alpha}$ 
to be the element such that 
\begin{align}\label{eq:flbb}
\iota_{\beta} (f_{\lambda}^{(b,\beta)})=P_\lambda^{(b)},
\end{align}
where $\iota_{\beta}$ is the isomorphism given in \eqref{eq:iota}.

(2)
For a complex number $\beta$, define an element of $\widehat{U}(\calH)$ by
\begin{align}\label{eq:E:split}
\widehat{E}_{\beta}
=\sqrt{2 \beta} \sum_{n>0}a_{-n} \calL_n
+\sum_{n>0}a_{-n}a_n\left(\beta-1-\sqrt{2\beta}a_0\right).
\end{align}
Here $\calL_n\in\widehat{U}(\calH)$ is the bosonized Virasoro generator 
\eqref{eq:FF}, and we have put the assumption
\begin{align}\label{eq:prp:split:1}
\rho=(\beta^{1/2}-\beta^{-1/2})/\sqrt{2}.
\end{align} 
\end{dfn}

\begin{prp}\label{prp:split}
For a generic complex number $\beta$ we have
\begin{align}
\label{eq:eigen:Ef}
&\widehat{E}_{\beta} f_{\lambda}^{(\beta^{-1},\beta)}
=\ep_{\lambda}(\beta)f_{\lambda}^{(\beta^{-1},\beta)},\\
\label{eq:eigen:E}
&\ep_{\lambda}(\beta)\seteq \sum_i (\lambda_i^2+\beta(1-2i)\lambda_i),
\end{align}
for any partition $\lambda$.
\end{prp}

The proof of this proposition is rather complicated,
since we should utilise Jack symmetric polynomials with finite variables.

\subsection{Jack symmetric polynomials}

Recall that in \S \ref{subsec:Fs} we denoted by $\Lambda_N$ 
the space of symmetric polynomials of $N$ variables,
and by $\Lambda_N^d$ its degree $d$  homogeneous subspace.
In order to denote $N$-variable symmetric polynomials, 
we put the superscript ``$(N)$" on the symbols 
for the infinite-variable symmetric functions.
For example,  we denote
by $p_\lambda^{(N)}(x)\seteq 
 p_{\lambda_1}^{(N)}(x) p_{\lambda_2}^{(N)}(x) \cdots$ 
the product of the power sum polynomials 
$p_k^{(N)}(x)\seteq \sum_{i=1}^N x_i^k$, 
and by $m_\lambda^{(N)}(x)$ 
the monomial symmetric polynomial.

Let us fix $N\in\bbZ_{\ge1}$ and an indeterminate
\footnote{In the literature this indeterminate is usually denoted by 
$\beta=\alpha^{-1}$, and we will also identify it with our $\beta$ given in 
\eqref{eq:prp:split:1} later.
But at this moment we don't use it to avoid confusion.} 
 $t$.
For a partition $\lambda$ with $\ell(\lambda)\le N$, 
the $N$-variable Jack symmetric polynomial 
$P_\lambda^{(N)}(x;t)$ is 
uniquely specified by the following two properties.
\begin{description}
\item[(i)] 
\begin{align*}
P_\lambda^{(N)}(x;t)
=m_\lambda^{(N)}(x)
+\sum_{\mu<\lambda}\widetilde{c}_{\lambda,\mu}(t)m_\mu^{(N)}(x),
\quad
\widetilde{c}_{\lambda,\mu}(t)\in\bbQ(t).
\end{align*}
\item[(ii)]
\begin{align}
\label{eq:eigeneq}
&H_t^{(N)} P_\lambda^{(N)}(x;t)
=\ep_{\lambda}^{(N)}(t) P_\lambda^{(N)}(x;t),
\\
\label{eq:HCS}
&\quad
H_{t}^{(N)}\seteq
 \sum_{i=1}^{N}\big(x_i\dfrac{\partial}{\partial x_i}\big)^2
 +t\sum_{1\le i<j\le N}\dfrac{x_i+x_j}{x_i-x_j}
  (x_i\dfrac{\partial}{\partial x_i}-x_j\dfrac{\partial}{\partial x_j}),
\\
&\quad
 \ep_{\lambda}^{(N)}(t)\seteq \sum_i (\lambda_i^2+t(N+1-2i)\lambda_i).
\label{eq:ep:N}
\end{align}
\end{description}
The differential operator \eqref{eq:HCS} is known 
to be equivalent to the Calogero-Sutherland Hamiltonian
(see \cite[\S 2]{AMOS:1995} for the detailed explanation.)
In (i) we used the dominance partial semi-ordering on the partitions.
If $N\ge d$, then $\{P_\lambda^{(N)}(x;t)\}_{\lambda\vdash d}$ is a basis 
of $\Lambda_{N,\bbQ(t)}^{d}$.

\begin{dfn}
For $M\ge N$, 
we denote the restriction map from $\Lambda_M$ to $\Lambda_N$ by 
\begin{align*}
\begin{array}{ c c c c}
\rho_{M,N}:  &\Lambda_M  &\to &\Lambda_N 
\\
       &\rotatebox{90}{$\in$} &  &\rotatebox{90}{$\in$} 
\\
       & f(x_1,\ldots,x_M)&\mapsto & f(x_1,\ldots,x_N,0,\ldots,0),
\end{array}
\end{align*}
and the induced restriction map from $\Lambda$ to $\Lambda_N$ by 
\begin{align*}
\begin{array}{ c c c c}
\rho_N: &\Lambda &\to &\Lambda_N.
\\
       &\rotatebox{90}{$\in$} &  &\rotatebox{90}{$\in$} 
\\
       & f(x_1,x_2,\ldots)&\mapsto & f(x_1,\ldots,x_N,0,\ldots).
\end{array}
\end{align*}
We denote the maps on the tensored spaces 
$\Lambda_{M,\bbC} \to \Lambda_{N,\bbC}$ and
$\Lambda_\bbC \to \Lambda_{N,\bbC}$ 
by the same symbols $\rho_{M,N}$ and $\rho_N$. 
\end{dfn}

\begin{fct}
For any $\lambda\in\calP$, every $N\in\bbZ_{\ge1}$ with $N\ge\ell(\lambda)$,
and any generic $t\in\bbC$ we have
\begin{align*}
\rho_{N} \big(P_{\lambda}^{(t^{-1})}\big) = P_\lambda^{(N)}(x;t).
\end{align*}
\end{fct}

\subsection{Split form of the Calogero-Sutherland Hamiltonian}

We recall the collective field method in the Calogero-Sutherland model 
following \cite[\S 3]{AMOS:1995}.
Recall that the Calogero-Sutherland Hamiltonian $H_t^{(N)}$ \eqref{eq:HCS}
acts on the space of symmetric polynomials $\Lambda_{N,\bbQ(t)}$.

\begin{fct}\label{fct:split}
(1)
Let $t,t'$ be non-zero complex numbers.
Define an element of $\widehat{U}(\calH)$ by 
\begin{align*}
\widehat{H}_{t,t'}^{(N)}\seteq
 \sum_{n,m>0}\Big(t' a_{-m-n}a_m a_n
                 &+\dfrac{t}{t'}a_{-m}a_{-n}a_{m+n}\Big)
\\
&+\sum_{n>0}\left(n(1-t)+N t\right) a_{-n}a_n.
\end{align*}
Then for any $v\in\calF_\alpha$ and every $N\in\bbZ_{\ge1}$ we have
\begin{align}\label{eq:fct:H}
\rho_{N}\circ\iota_{t} (\widehat{H}_{t,t'}^{(N)} v) 
=H_{t}^{(N)}\big( \rho_{N}\circ \iota_{t}(v)\big).
\end{align}

(2)
Under the relation
\[\rho=\big(t^{1/2}-t^{-1/2}\big)/\sqrt{2}\]
we have
\begin{align}\label{eq:fct:split:2}
\widehat{H}_{t,\sqrt{t/2}}^{(N)}
=\sqrt{2 t} \sum_{n>0}a_{-n} \calL_n
+\sum_{n>0}a_{-n}a_n\left(N t+t-1-\sqrt{2 t} a_0\right).
\end{align}
Here $\calL_n\in\widehat{U}(\calH)$ is the bosonized Virasoro generator 
\eqref{eq:FF}.
\end{fct}

\begin{proof}
These are well-known results 
(for example, see \cite[Prop. 4.47]{S:2003}, \cite{AMOS:1995} 
 and the references therein). 
We only show the sketch of the proof. 

As for (1),
note that $\{a_{-\lambda}\ket{\alpha}_{\calF}\mid \lambda\in\calP\}$ 
is a basis of $\calF$.
So it is enough to show \eqref{eq:fct:H} for each $\lambda$.
One can calculate the left hand side using the commutation 
relation of $\calH$ only.
On the right hand side, 
one may use $\iota_t(a_{-\lambda}\ket{\alpha}_{\calF})\propto p_\lambda$,
and calculate it using the expression \eqref{eq:HCS}.

(2) is proved by direct calculation.
\end{proof}

\begin{rmk}
The form \eqref{eq:fct:split:2} is called the split expression in \cite[\S 1]{SKAO:1996}.
\end{rmk}

\subsection{Proof of Proposition \ref{prp:split}}

By Fact \ref{fct:split} we have 
the left commuting diagram in \eqref{eq:comm:diag}.
Note that we set the parameters $t$ and $t'$ 
in $\widehat{H}_{t,t'}^{(N)}$ 
to be $\beta$ and $\sqrt{\beta/2}$, 
so that we may use Fact \ref{fct:split} (2).
In the right diagram of \eqref{eq:comm:diag} we show how 
the element $f_\lambda^{(b,\beta)}\in\calF_\alpha$ 
given in \eqref{eq:flbb} behaves under the maps indicated 
in the left diagram.
Here we set the parameter $b$ to be $\beta^{-1}$ so that 
$\iota_\beta f_\lambda^{(\beta^{-1},\beta)}
=P_\lambda^{(\beta^{-1})} \in\Lambda_\bbC$
and 
$\rho_N\circ\iota_\beta f_\lambda^{(\beta^{-1},\beta)}
=P_\lambda^{(N)}(x;\beta) \in\Lambda_{N,\bbC}$.
At the bottom line we used 
the eigen-equation of Jack symmetric polynomial \eqref{eq:eigeneq}.

\begin{align}
\xymatrix{
  \calF_\alpha  
   \ar[r]^{\widehat{H}_{\beta,\sqrt{\beta/2}}^{(N)}}
   \ar[d]_{\iota_\beta}^{\rotatebox{90}{$\sim$}}
& \calF_\alpha 
   \ar[d]^{\iota_\beta}_{\rotatebox{90}{$\sim$}}  
& f_\lambda^{(\beta^{-1},\beta)} 
   \ar@{|->}[d]
   \ar@{|->}[r]
& \widehat{H}_{\beta,\sqrt{\beta/2}}^{(N)} 
  \big(f_\lambda^{(\beta^{-1},\beta)} \big)
   \ar@{|->}[d]
\\
  \Lambda_\bbC  
   \ar[d]_{\rho_N}
   \ar@{}[r]|{\circlearrowright}
& \Lambda_\bbC 
   \ar[d]^{\rho_N}
& P_\lambda^{(\beta^{-1})} 
   \ar@{|->}[d]
& \iota_\beta\circ\widehat{H}_{\beta,\sqrt{\beta/2}}^{(N)} 
  \big(f_\lambda^{(\beta^{-1},\beta)}\big)
   \ar@{|->}[d]
\\ 
  \Lambda_{N,\bbC}
   \ar[r]_{H_{\beta}^{(N)}}
& \Lambda_{N,\bbC}
& P_\lambda^{(N)}(x;\beta) 
   \ar@{|->}[r]
& P_\lambda^{(N)}(x;\beta)\cdot \ep_\lambda^{(N)}(\beta)
}
\label{eq:comm:diag}
\end{align}

Since this diagram holds for every $N$ with $N\ge\ell(\lambda)$, we have
\begin{align*}
\widehat{H}_{\beta,\sqrt{\beta/2}}^{(N)} f_\lambda^{(\beta^{-1},\beta)}
=\ep_\lambda^{(N)}(\beta) f_\lambda^{(\beta^{-1},\beta)}.
\end{align*}
Therefore we have
\begin{align*}
&
\big[\sqrt{2 t} \sum_{n>0}a_{-n} \calL_n
+\sum_{n>0}a_{-n}a_n\left(N t+t-1-\sqrt{2 t} a_0\right)\big]
 f_\lambda^{(\beta^{-1},\beta)}
\\
&=f_\lambda^{(\beta^{-1},\beta)}\cdot
 \sum_i (\lambda_i^2+t(N+1-2i)\lambda_i).
\end{align*}
We can subtract $N$-dependent terms from both sides.
The result is nothing but the desired statement 
of Proposition \ref{prp:split}.

\section{Whittaker vectors}\label{sec:whit}

Recall the notion of the Whittaker vector for a finite dimensional 
Lie algebra $\frkg$ given in \cite{K:1978}. 
Let $\frkn$ be a maximal nilpotent Lie subalgebra of $\frkg$ 
and $\eta:\frkn\to\bbC$ be a homomorphism.
Let $V$ be any $U(\frkg)$-module.
Then a vector $w\in V$ is called a Whittaker vector with respect to $\eta$ if 
$x w=\eta(x)w$ for all $x\in\frkn$.

We shall discuss an analogue of this Whittaker vector 
in the Virasoro algebra $\Vir_c$.
In the triangular decomposition 
$\Vir_{c}=\Vir_{c,+}\oplus\Vir_{c,0}\oplus\Vir_{c,-}$,
the elements $L_1,L_2\in\Vir_{c,+}$ generate $\Vir_{c,+}$.
Thus if we take $\Vir_{c,+}$ as the $\eta$ in the above definition,
what we should consider is a homomorphism $\eta:\Vir_{c,+}\to\bbC$,
which is determined by $\eta(L_1)$ and $\eta(L_2)$.

In \cite{OW:2009}, a characterisation of Whittaker vectors 
in general $U(\Vir)$-modules are given
under the assumption that $\eta$ is non-degenerate, 
i.e. $\eta(L_1)\neq0$ and $\eta(L_2)\neq0$.

In this section we shall express Whittaker vectors 
in the Verma module $M_{h}$ using Jack symmetric functions.
Before starting the general treatment, 
we first investigate a degenerate version of the Whittaker vector,
i.e. we assume $\eta(L_2)=0$.
We will call this vector by Gaiotto state of Virasoro algebra,
although the paper \cite{G:2009} treated both degenerate and 
non-degenerate Whittaker vectors.

\subsection{Gaiotto state via Jack polynomials}\label{subsec:coh}

\begin{dfn}
Fix a non-zero complex number $a$.
Let $v_G$ be a non-zero element of the Verma module $M_{h}$ satisfying
\begin{align*}
L_1 v_G=a v_G,\quad
L_n v_G=0\ (n\in\bbZ_{\ge 2}).
\end{align*}
We call such an element $v_G$ a Gaiotto state of $M_{h}$.
\end{dfn}

\begin{fct}\label{fct:G:ue}
Assume that $c$ and $h$ are generic.
Then $v_G$ exists uniquely up to constant multiplication.
\end{fct}
\begin{proof}
This statement is shown in \cite{OW:2009}.
\end{proof}

\begin{lem}
Decompose a Gaiotto state $v_G$ in the way \eqref{eq:L0weightdecomp} as 
\begin{align*}
v_G=\sum_{n\in\bbZ_{\ge 0}}a^n v_{G,n},\quad
v_{G,n}\in M_{h,n}.
\end{align*}
Then we have
\begin{align}\label{eq:cond:g}
v_{G,n}=L_1 v_{G,n+1} \quad (n\in\bbZ_{\ge 0}).
\end{align}
\end{lem}

\begin{proof}
This follows from the commutation relation $[L_0,L_1]=-L_1$.
\end{proof}

Now consider the bosonized Gaiotto state 
\begin{align*}
w_{G,n}\seteq \psi(v_{G,n})\in\calF_{\alpha,n}
\end{align*}
where $\psi: M_h \to \calF_{\alpha}$ 
is the Feigin-Fuchs bosonization \eqref{eq:FF:psi}
and $\calF_{\alpha,n}$ is the $a_0$-weight space \eqref{eq:F:deg}.
At this moment the Heisenberg parameters $\rho,\alpha$ 
are related to the Virasoro parameters $c,h$ by the relations
\begin{align*}
c=1-12\rho^2,\quad
h=\dfrac{1}{12}\alpha(\alpha-2\rho).
\end{align*}
From the condition \eqref{eq:cond:g} we have
\begin{align}\label{eq:cond:wg}
\calL_1 w_{G,n+1}\in \calF_{h,n},\quad w_{G,n}=\calL_1 w_{G,n+1}.
\end{align}

Next we map this bosonized state $w_{G,n}$ into a symmetric function 
by the isomorphism $\iota_\beta:\calF_{\alpha}\to\Lambda_\bbC$ \eqref{eq:iota}:
\begin{align*}
\iota_\beta(w_{G,n})=\iota_\beta\circ\psi(v_{G,n})\in\Lambda_\bbC^n.
\end{align*}
Here $\Lambda_\bbC^n$ is the space of degree $n$ symmetric functions.
We take the parameter $\beta$ so that 
the Heisenberg parameter $\rho$ is expressed by 
\begin{align*}
\rho=(\beta^{1/2}-\beta^{-1/2})/\sqrt{2}.
\end{align*}
Recall also that the family of Jack symmetric functions 
$\{P_\lambda^{(\beta^{-1})} \mid \lambda\vdash n\}$ 
is a basis of $\Lambda_{\bbC}^n$ for a generic $\beta\in\bbC$.
Thus we can expand $\iota_\beta(w_{G,n})\in\Lambda_\bbC^n$ 
by $P_\lambda^{(\beta^{-1})}$'s.
Let us express this expansion as:
\begin{align}\label{eq:expand}
\iota_{\beta}(w_{G,n})=\iota_{\beta}\circ\psi(v_{G,n})
=\sum_{\lambda\vdash n}c_\lambda(\alpha,\beta)P_{\lambda}^{(\beta^{-1})},\quad
c_\lambda(\alpha,\beta)\in\bbC.
\end{align}
Note that this expansion is equivalent to 
\begin{align}\label{eq:expand:Fock}
w_{G,n}
=\sum_{\lambda\vdash n}c_\lambda(\alpha,\beta)f_{\lambda}^{(\beta^{-1},\beta)}
\in \calF_\alpha
\end{align}
by \eqref{eq:flbb}.
Now the correspondence of the parameters becomes: 
\begin{align}\label{eq:param}
c=1-12\rho^2,\quad
h=\dfrac{1}{12}\alpha(\alpha-2\rho_0),\quad
\rho=\dfrac{\beta^{1/2}-\beta^{-1/2}}{\sqrt{2}}.
\end{align}

The main result of this paper is 
\begin{thm}\label{thm:G}
Assume that $c$ and $h$ are generic.
(Then $v_G$ exists  uniquely up to constant multiplication by Fact \ref{fct:G:ue}.)
If $c_{\emptyset}(\alpha,\beta)$ is set to be $1$ 
in the expansion \eqref{eq:expand}, 
then the other coefficients are given by
\begin{align}\label{eq:thm:G}
\begin{split}
c_\lambda(\alpha,\beta)
=\prod_{(i,j)\in\lambda}
 &\dfrac{1}{\lambda_i-j+1+\beta(\lambda_j'-i)}
\\
&\times
 \prod_{(i,j)\in\lambda\setminus\{(1,1)\}}
 \dfrac{1}{(j+1)\beta+\sqrt{2} \beta^{1/2}\alpha -(i+1)}.
\end{split}
\end{align}
Here we used the notation $(i,j)\in\lambda$ as explained in \S \ref{subsec:partition}.
\end{thm}

\subsection{Proof of Theorem \ref{thm:G}}

Before starting the proof, we need to prepare the following 
Proposition \ref{prp:rec}.
Recall the Pieri formula of Jack symmetric function.
We only need the case of ``adding one box", that is, 
the case of multiplying the degree one power sum function $p_1=x_1+x_2+\cdots$.

\begin{dfn}\label{dfn:<_k}
For partitions $\mu$ and $\lambda$,
we denote $\mu<_k\lambda$ if $|\mu|=|\lambda|-k$ and $\mu\subset\lambda$.
\end{dfn}

\begin{fct}[{\cite[p.340 VI (6.24), p.379 VI (10.10)]{M:1995:book}}]
We have
\begin{align}
\label{eq:pieri}
& p_1 P_\mu^{(\beta^{-1})} 
=\sum_{\lambda>_1\mu}  \psi_{\lambda/\mu}'(\beta) P_\lambda^{(\beta^{-1})},
\\
&\label{eq:pieri:coeff}
\psi_{\lambda/\mu}'(\beta)\seteq
\prod_{i=1}^{I-1}
\dfrac{\lambda_i-\lambda_I+\beta(I-i+1)}{\lambda_i-\lambda_I+1+\beta(I-i)}
\dfrac{\lambda_i-\lambda_I+1+\beta(I-i-1)}{\lambda_i-\lambda_I+\beta(I-i)}.
\end{align}
In the expression in \eqref{eq:pieri:coeff} 
the partitions $\lambda$ and $\mu$ are related by $\lambda_I=\mu_I+1$ and 
$\lambda_i=\mu_i$ for $i\neq I$.
\end{fct}

\begin{prp}\label{prp:rec}
$c_\lambda(\alpha,\beta)$ satisfies the next recursion relation.
\begin{align}\label{eq:prp:rec}
\left(\ep_{\lambda}(\beta)+|\lambda|(1+\sqrt{2\beta}\alpha-\beta)\right)
 c_\lambda(\alpha,\beta)
=\beta \sum_{\mu<_1\lambda}\psi_{\lambda/\mu}'(\beta)c_\mu(\alpha,\beta).
\end{align}
Here the function $\ep_{\lambda}(\beta)$ is given in \eqref{eq:eigen:E}.
\end{prp}

\begin{proof}
We will calculate $\widehat{E}_\beta w_{G,n}\in\calF_\alpha$ in two ways.
By comparing both expression we obtain the recursion relation.

First, by the definition of $\widehat{E}_\beta$ given in \eqref{eq:E:split} 
and by the condition \eqref{eq:cond:wg} of $v_{G,n}$ we have
\begin{align*}
\widehat{E}_\beta w_{G,n}
&=\Big[\sqrt{2 \beta} \sum_{m\ge 1}a_{-m}\calL_m 
      +\sum_{m\ge 1} a_{-m}a_m(\beta-1-\sqrt{2 \beta}a_0)
  \Big]w_{G,n}\\
&=\Big[\sqrt{2 \beta} a_{-1}\calL_1
      +\sum_{m\ge 1} a_{-m}a_m(\beta-1-\sqrt{2 \beta}a_0)
  \Big]w_{G,n}\\
&=\sqrt{2\beta}a_{-1}w_{G,n-1}+n(\beta-1-\sqrt{2\beta}\alpha)
  w_{G,n} \in \calF_\alpha.
\end{align*}
Now applying the isomorphism $\iota_\beta:\calF_\alpha\to\Lambda_\bbC$ 
on both sides and 
substituting $w_{G,n}$ and $w_{G,n-1}$ by their expansions \eqref{eq:expand},
we have
\begin{align*}
\iota_\beta(\widehat{E}_\beta w_{G,n})
= \beta p_{1}
&\sum_{\mu\vdash n-1} c_{\mu}(\alpha,\beta)P_{\mu}^{(\beta^{-1})}
\\
&+n(\beta-1-\sqrt{2\beta}\alpha)
  \sum_{\lambda\vdash n}c_\lambda(\alpha,\beta) P_{\lambda}^{(\beta^{-1})}
 \in \Lambda_\bbC.
\end{align*}
Using the Pieri formula \eqref{eq:pieri} in the first term, we have
\begin{align}\label{eq:prp:rec:1}
\begin{split}
\iota_\beta(\widehat{E}_\beta w_{G,n})
=&\beta \sum_{\mu\vdash n-1}c_{\mu}(\alpha,\beta)
    \sum_{\lambda>_1\mu}\psi'_{\lambda/\mu}(\beta) P_{\lambda}^{(\beta^{-1})}\\
 &+n(\beta-1-\sqrt{2\beta}\alpha)\sum_{\lambda\vdash n}
  c_\lambda(\alpha,\beta) P_\lambda^{(\beta^{-1})}.
\end{split}
\end{align}

Next, by \eqref{eq:expand:Fock} and by \eqref{eq:eigen:Ef} we have 
\begin{align*}
\widehat{E}_\beta w_{G,n}
=\widehat{E}_\beta \sum_{\lambda\vdash n}c_\lambda(\alpha,\beta)
 f_{\lambda}^{(\beta^{-1},\beta)}
=\sum_{\lambda\vdash n}
 c_\lambda(\alpha,\beta) \ep_{\lambda}(\beta)
 f_{\lambda}^{(\beta^{-1},\beta)}\in\calF_\alpha.
\end{align*}
Therefore we have 
\begin{align}\label{eq:prp:rec:2}
\iota_\beta(\widehat{E}_\beta w_{G,n})
=\sum_{\lambda\vdash n}
 c_\lambda(\alpha,\beta) \ep_{\lambda}(\beta) P_\lambda^{(\beta^{-1})}
 \in\Lambda_\bbC.
\end{align}
Then comparing \eqref{eq:prp:rec:1} and \eqref{eq:prp:rec:2} we have
\begin{align*}
&\sum_{\lambda\vdash n}
 \big(\ep_{\lambda}(\beta)+n(1+\sqrt{2\beta}\alpha-\beta)\big)
 c_\lambda(\alpha,\beta)
 P_\lambda^{(\beta^{-1})}
\\
&=\beta\sum_{\mu\vdash n-1}c_\mu(\alpha,\beta)
  \sum_{\lambda>_1\mu}\psi'_{\lambda/\mu}(\beta) P_\lambda^{(\beta^{-1})}
  \in\Lambda_{\bbC}^n.
\end{align*}
Since $\{P_\lambda^{(\beta)} \mid \lambda\vdash n\}$ 
is a basis of $\Lambda_{\bbC}^n$, 
we have 
\begin{align*}
\big(\ep_{\lambda}(\beta)+n(1+\sqrt{2\beta}\alpha-\beta)\big)
c_\lambda(\alpha,\beta)
=\beta \sum_{\mu<_1\lambda}c_\mu(\alpha,\beta)\psi'_{\lambda/\mu}(\beta).
\end{align*}
\end{proof}

\begin{proof}[Proof of Theorem \ref{thm:G}]
The recursion relation \eqref{eq:prp:rec} of Propositions \ref{prp:rec} 
determines $c_\lambda(\alpha,\beta)$ uniquely 
if we set the value of $c_{\emptyset}(\alpha,\beta)$.
Since the existence and uniqueness of $v_G$ is known by Fact \ref{fct:G:ue},
we only have to show that the ansatz \eqref{eq:thm:G} 
satisfies \eqref{eq:prp:rec}.

For partitions $\lambda$ and $\mu$ which are related 
by $\lambda_I=\mu_I+1$ and $\lambda_i=\mu_i$ for $i\neq I$, 
we have the following two formulas:
\begin{align*}
&\Big[\prod_{(i,k)\in\mu}\dfrac{1}{\lambda_i-k+1+\beta(\lambda_k'-i)}\Big]
 \Big/
 \Big[\prod_{(i,k)\in\lambda}\dfrac{1}{\lambda_i-k+1+\beta(\lambda_k'-i)}\Big]
\\
&
=\prod_{i=1}^{I-1}
 \dfrac{\lambda_i-\lambda_I+1+\beta(I-i)}
       {\lambda_i-\lambda_I+1+\beta(I-1-i)}
 \times
 \prod_{i=1}^{\lambda_I-1}
 \dfrac{\lambda_I-i+1+\beta(\lambda_i'-I)}{\lambda_I-i+\beta(\lambda_i'-I)},
\\
&\Big[\prod_{\substack{(i,k)\in\mu \\ \mu\neq (1,1) }}
      \dfrac{\beta}{(k+1)+\sqrt{2\beta}\alpha -(i+1)\beta}\Big]
\Big/
 \Big[\prod_{\substack{(i,k)\in\lambda \\ \mu\neq (1,1)}}
      \dfrac{\beta}{(k+1)+\sqrt{2\beta}\alpha -(i+1)\beta}\Big]
\\
&
=\dfrac{(\lambda_I+1)+\sqrt{2\beta}\alpha -(I+1)\beta}{\beta}.
\end{align*}
Substituting the $c_\mu(\alpha,\beta)$ in the right hand side of 
\eqref{eq:prp:rec} by the ansatz \eqref{eq:thm:G} 
and using the above two equations, we have
\begin{align}
\label{eq:prf:thm:rhs}
\begin{split}
&\text{RHS of \eqref{eq:prp:rec}}
\\
&=
 \sum_{(I,\lambda_I)\in C(\lambda)} 
 \prod_{i=1}^{I-1}
 \dfrac{\lambda_i-\lambda_I+\beta(I-i+1)}
       {\lambda_i-\lambda_I+\beta(I-i)}
 \times
 \prod_{i=1}^{\lambda_I-1}
 \dfrac{\lambda_I-i+1+\beta(\lambda_i'-I)}{\lambda_I-i+\beta(\lambda_i'-I)}
\\
&\hskip 6em
 \times
 \left((\lambda_I+1)+\sqrt{2\beta}\alpha -(I+1)\beta\right)
 c_\lambda(\alpha,\beta),
\end{split}
\end{align}
where $C(\lambda)$ is the set of boxes $\square$ 
in the Young diagram of $\lambda$ 
such that $\lambda\setminus \{\square\}$ is also a partition.
In particular, if $\square=(I,\lambda_I)\in C(\lambda)$, then 
$\mu\seteq\lambda\setminus \{\square\}$ is the partition satisfying 
$\mu_I=\lambda_I-1$ and $\mu_i=\lambda_i$ for $i\neq I$,
recovering the previous description.

As for the left hand side of \eqref{eq:prp:rec}, 
we have by \eqref{eq:eigen:E}:
\begin{align}
\label{eq:prf:thm:lhs}
\ep_\lambda(\beta)+|\lambda|(1+\sqrt{2\beta}\alpha-\beta)
=|\lambda|(1+\sqrt{2\beta}\alpha)+\sum_i (\lambda_i^2-2i\lambda_i\beta).
\end{align}

Thus by \eqref{eq:prf:thm:rhs} and \eqref{eq:prf:thm:rhs},
the equation \eqref{eq:prp:rec} under the substitution 
\eqref{eq:thm:G} is equivalent to the next one:
\begin{align*}
&|\lambda|(1+\sqrt{2\beta}\alpha)+\sum_i (\lambda_i^2-2i\lambda_i\beta)
\\
&=\sum_{(I,\lambda_I)\in C(\lambda)} 
 \prod_{i=1}^{I-1}
 \dfrac{\lambda_i-\lambda_I+\beta(I-i+1)}
       {\lambda_i-\lambda_I+\beta(I-i)}
 \times
 \prod_{i=1}^{\lambda_I-1}
 \dfrac{\lambda_I-i+1+\beta(\lambda_i'-I)}{\lambda_I-i+\beta(\lambda_i'-I)}
\\
&\phantom{=c_\lambda(\alpha,\beta) \sum_{\mu<_1 \lambda}\times }
 \times
 \left(1+\sqrt{2\beta}\alpha+\lambda_I-(I+1)\beta\right).
\end{align*}
This is verified by Propositions \ref{prp:young:1} and \ref{prp:young:2} 
which will be shown in the next \S \ref{sec:id}.
\end{proof}

\subsection{Non-degenerate Whittaker vector via Jack polynomials}\label{subsec:wit}

\begin{dfn}
Fix non-zero complex numbers $a$ and $b$.
Let $v_W$ be an element of the Verma module $M_{h}$ satisfying
\begin{align*}
L_1 v_W=a v_W,\quad
L_2 v_W=b v_W,\quad
L_n v_W=0\ (n\in\bbZ_{\ge 3}).
\end{align*}
We call such an element $v_W$ by (non-degenerate) Whittaker vector of $M_{h}$.
\end{dfn}

\begin{fct}
Assume that $c$ and $h$ are generic complex numbers. 
Then $v_W$ exists uniquely up to scalar multiplication.
\end{fct}

\begin{proof}
This is shown in \cite{OW:2009}.
\end{proof}

\begin{lem}
Let us decompose $v_W$ as
\begin{align*}
v_W=\sum_{n\in\bbZ_{\ge 0}}a^n v_{W,n},\quad
v_{W,n}\in M_{h,n}.
\end{align*}
Then we have 
\begin{align*}
&L_1 v_{W,n+1}=v_{W,n},\quad
 L_2 v_{W,n+2}=a^{-2} b\cdot v_{W,n}.
\end{align*}
\end{lem} 
\begin{proof}
This follows from the commutation relations $[L_0,L_1]=-L_1$ and 
$[L_0,L_2]=-2L_2$.
\end{proof}

Now we expand the bosonized Whittaker vector
\begin{align*}
w_{W,n}\seteq \psi(v_{W,n})\in\calF_{\alpha,n}
\end{align*} 
by $f_{\lambda}^{(\beta^{-1},\beta)}$'s \eqref{eq:flbb} 
and express it as
\begin{align*}
w_{W,n}
=\sum_{\lambda\vdash n}d_\lambda(\alpha,\beta)f_{\lambda}^{(\beta^{-1},\beta)},
\quad d_\lambda(\alpha,\beta)\in\bbC.
\end{align*}

\begin{prp}\label{prp:rec2}
Using the notation $\lambda >_k \mu$ given in Definition \ref{dfn:<_k},
we have the next recursion relation for $d_\lambda(\alpha,\beta)$:
\begin{align}\label{eq:prp:rec2}
\begin{split}
&(\ep_{\lambda}(\beta)+|\lambda|(1+\sqrt{2\beta}\alpha-\beta))
c_\lambda(\alpha,\beta)
\\
&=\beta \sum_{\nu<_2\lambda}
   d_\nu(\alpha,\beta){\psi'}_{\lambda/\nu}^{(2)}(\beta) 
 +\beta \sum_{\mu<_1\lambda}d_\mu(\alpha,\beta)\psi'_{\lambda/\mu}(\beta).
\end{split}
\end{align}
${\psi'}_{\lambda/\nu}^{(2)}(\beta)$ is the coefficient 
in the next Pieri formula:
\begin{align*}
p_2 P_{\nu}^{(\beta^{-1})}
=\sum_{\lambda} {\psi'}_{\lambda/\nu}^{(2)}(\beta) P_{\lambda}^{(\beta^{-1})}.
\end{align*}
\end{prp}

\begin{proof}
Similar as the proof of Proposition \ref{prp:rec}
\end{proof}

\begin{rmk}
The author doesn't know whether $d_\lambda$ has a good explicit formula,
although $c_\lambda$ has the factored formula \eqref{eq:thm:G}.
\end{rmk}

\section{Combinatorial identities of rational functions}\label{sec:id}

\begin{prp}\label{prp:young:1}
For a partition $\lambda$, 
let $C(\lambda)$ be the set of boxes $\square$ of $\lambda$ 
such that $\lambda\setminus \{\square\}$ is also a partition.
Then
\begin{align}\label{eq:young:1}
 \sum_{(I,\lambda_I)\in C(\lambda)} 
 \prod_{i=1}^{I-1}
 \dfrac{\lambda_i-\lambda_I+\beta(I-i+1)}
       {\lambda_i-\lambda_I+\beta(I-i)}
 \times
 \prod_{i=1}^{\lambda_I-1}
 \dfrac{\lambda_I-i+1+\beta(\lambda_i'-I)}{\lambda_I-i+\beta(\lambda_i'-I)}
=|\lambda|.
\end{align}
\end{prp}

\begin{proof}
Let $\lambda$ be the partition such that 
\begin{align}\label{eq:young:1:lambda}
\lambda=(
 \stackrel{j_1}{\overbrace{n_1,\ldots,n_1}},
 \stackrel{j_2}{\overbrace{n_2,\ldots,n_2}},\ldots,
 \stackrel{j_l}{\overbrace{n_l,\ldots,n_l}}).
\end{align} 
Then we have
\begin{align*}\label{eq:young:1:C}
C(\lambda)=\{(m_1,n_1),(m_2,n_2),\ldots,(m_l,n_l)\}
\end{align*} 
with $m_k\seteq j_1+\cdots+j_k$ ($k=1,\ldots,l)$,
where we used the coordinate $(i,j)$ of Young diagram associated to $\lambda$ 
as explained in \S \ref{subsec:partition}.

Let us choose an element $\square=(m_k,n_k)$ of $C(\lambda)$, 
and calculate the corresponding factor in \eqref{eq:young:1}.
The first product reads
\begin{align*}
&\prod_{1\le i\le m_1}\dfrac{(n_1-n_k)+\beta(m_k-i+1)}{(n_1-n_k)+\beta(m_k-i)}
\prod_{m_1+1\le j\le m_2}
\dfrac{(n_2-n_k)+\beta(m_k-i+1)}{(n_2-n_k)+\beta(m_k-i)}
\\
&\times\cdots\times
\prod_{m_{k-1}+1\le i\le m_{k}-1}\dfrac{m_k-i+1}{m_k-i}
\\
&=\prod_{i=1}^{k-1}
 \dfrac{(n_i-n_k)+\beta(m_k-m_{i-1})}{(n_i-n_k)+\beta(m_k-m_i)}
 \times(m_k-m_{k-1}).
\end{align*}
Here we used the notation $m_0\seteq0$.
The second product reads
\begin{align*}
&\prod_{1\le j\le n_l}\dfrac{(n_k-j+1)+\beta(m_l-n_k)}{(n_k-j)+\beta(m_l-n_k)}
\prod_{n_l+1\le j\le n_{l-1}}
\dfrac{(n_k-j+1)+\beta(m_{l-1}-n_k)}{(n_k-j)+\beta(m_{l-1}-n_k)}
\\
&\times\cdots\times
\prod_{n_{k+1}+1\le j\le n_{k}-1}
\dfrac{(n_k-j+1)}{(n_k-j)}
\\
&
=(n_k-n_{k+1})\times\prod_{j=k+1}^{l}
 \dfrac{(n_k-n_{j+1})+\beta(m_j-m_k)}{(n_k-n_j)+\beta(m_j-m_k)}.
\end{align*}
Here we used the notation $n_{l+1}\seteq0$.

Now let us define
\begin{align*}
&F_1(\{m_k\},\{n_k\},\beta)
\seteq \sum_{k=1}^l F_{1,k}(\{m_k\},\{n_k\},\beta),
\\
&F_{1,k}(\{m_k\},\{n_k\},\beta)
\seteq  (m_k-m_{k-1})(n_k-n_{k+1})
\\
&\times
 \prod_{i=1}^{k-1}
 \dfrac{(n_i-n_k)+\beta(m_k-m_{i-1})}{(n_i-n_k)+\beta(m_k-m_i)}
 \prod_{j=k+1}^{l}
 \dfrac{(n_k-n_{j+1})+\beta(m_j-m_k)}{(n_k-n_j)+\beta(m_j-m_k)}.
\end{align*}
Then for the proof of \eqref{eq:young:1}
it is enough to show that $F_1$ is equal to $|\lambda|$ 
if $\{m_k\}$ and $\{n_k\}$ correspond to $\lambda$ 
as in \eqref{eq:young:1:lambda} and \eqref{eq:young:1:C}.

Hereafter we consider $F_1$ as a rational function of the valuables 
$\{m_k\}$, $\{n_k\}$ and $\beta$.
As a rational function of $\beta$, $F_1$ has the apparent poles at 
$\beta_{j,k}\seteq -(n_j-n_k)/(m_k-m_j)$ ($j=1,2,\ldots,k-1,k+1,\ldots,l$).
We may assume that these apparent poles are mutually different 
so that all the poles are at most single.
Then the residue at $\beta=\beta_{j,k}$ comes from 
the factors $F_{1,j}$ and $F_{1,k}$.
Now we may assume $j<k$. Then the direct computation yields
\begin{align}
\Res_{\beta=\beta_{j,k}}F_{1,j}
=&\dfrac{(m_j-m_{j-1})(n_j-n_{j+1})(n_k-n_{k+1})}{(m_j-m_k)}
\\
\label{eq:prp:yng:1:j:1}
&\times
 \prod_{i=1}^{j-1}
 \dfrac{(n_i-n_j)(m_k-m_j)-(n_j-n_k)(m_j-m_{i-1})}
       {(n_i-n_j)(m_k-m_j)-(n_j-n_k)(m_j-m_i)}
\\
\label{eq:prp:yng:1:j:2}
&\times
 \prod_{i=j+1}^{k-1}
 \dfrac{(n_j-n_{i+1})(m_k-m_j)-(n_j-n_k)(m_i-m_j)}
       {(n_j-n_i)(m_k-m_j)-(n_j-n_k)(m_i-m_j)}
\\
\label{eq:prp:yng:1:j:3}
&\times
 \prod_{i=k+1}^{l}
 \dfrac{(n_j-n_{i+1})(m_k-m_j)-(n_j-n_k)(m_i-m_j)}
       {(n_j-n_i)(m_k-m_j)-(n_j-n_k)(m_i-m_j)},
\end{align}
and
\begin{align}
\Res_{\beta=\beta_{j,k}}F_{1,k}
=&\dfrac{(m_k-m_{k-1})(n_k-n_{k+1})(n_j-n_k)(m_{j-1}-m_j)}{(m_k-m_j)^2}
\\
\label{eq:prp:yng:1:k:1}
&\times
 \prod_{i=1}^{j-1}
 \dfrac{(n_i-n_k)(m_k-m_j)-(n_j-n_k)(m_k-m_{i-1})}
       {(n_i-n_k)(m_k-m_j)-(n_j-n_k)(m_k-m_i)}
\\
\label{eq:prp:yng:1:k:2}
&\times
 \prod_{i=j+1}^{k-2}
 \dfrac{(n_i-n_k)(m_k-m_j)-(n_j-n_k)(m_k-m_{i-1})}
       {(n_i-n_k)(m_k-m_j)-(n_j-n_k)(m_k-m_i)}
\\
\label{eq:prp:yng:1:k:3}
&\times
 \prod_{i=k+1}^{l}
 \dfrac{(n_k-n_{i+1})(m_k-m_j)-(n_j-n_k)(m_i-m_k)}
       {(n_k-n_i)(m_k-m_j)-(n_j-n_k)(m_i-m_k)}.
\end{align}
Using the identity $(a-b)(x-y)-(c-b)(x-z)=(a-c)(x-y)-(c-b)(y-z)$,
one finds that 
the factors \eqref{eq:prp:yng:1:j:1} and \eqref{eq:prp:yng:1:k:1} are equal.
Similarly \eqref{eq:prp:yng:1:j:3} and \eqref{eq:prp:yng:1:k:3} are equal.
Shifting the index $i$ in \eqref{eq:prp:yng:1:j:2} 
and using the above identity, 
one also finds that 
\begin{align*}
\dfrac{\eqref{eq:prp:yng:1:j:2}}{\eqref{eq:prp:yng:1:k:2}}
=-\dfrac{(n_j-n_k)(m_k-m_{k-1})}
        {(n_{j+1}-n_j)(m_k-m_j)}.
\end{align*}
Thus we have
\begin{align*}
\dfrac{\Res_{\beta=\beta_{j,k}}F_{1,j}}{\Res_{\beta=\beta_{j,k}}F_{1,k}}
=-\dfrac{(n_j-n_k)(m_k-m_{k-1})}{(n_{j+1}-n_j)(m_k-m_j)}
  \dfrac{\eqref{eq:prp:yng:1:j:1}}{\eqref{eq:prp:yng:1:k:1}}
=-1.
\end{align*}
Therefore we have
\begin{align*}
\Res_{\beta=\beta_{j,k}}F_1(\{m_k\},\{n_k\},\beta)=0,
\end{align*}
so that $F_1$ is a polynomial of $\beta$.

Then from the behaviour $F_1$ in the limit $\beta\to \infty$,
we find that $F_1$ is a constant as a function of $\beta$.
This constant can be calculated by setting $\beta=0$, and the result is
\begin{align*}
F_1(\{m_k\},\{n_k\},\beta)
=\sum_{k=1}^l  (m_k-m_{k-1})n_k.
\end{align*}
It equals to $|\lambda|$ if $\lambda$ is given by \eqref{eq:young:1:lambda}
and $m_k=j_1+\cdots+j_k$.
This is the desired consequence.
\end{proof}

\begin{prp}\label{prp:young:2}
Using the same notation as in Proposition \ref{prp:young:1}, we have
\begin{align}\label{eq:prp:young:2}
\begin{split}
&\sum_{(I,\lambda_I)\in C(\lambda)} 
 \prod_{i=1}^{I-1}
 \dfrac{\lambda_i-\lambda_I+\beta(I-i+1)}
       {\lambda_i-\lambda_I+\beta(I-i)}
 \times
 \prod_{i=1}^{\lambda_I-1}
 \dfrac{\lambda_I-i+1+\beta(\lambda_i'-I)}{\lambda_I-i+\beta(\lambda_i'-I)}
\\
&\phantom{=c_\lambda(\alpha,\beta) \sum_{\mu<_1 \lambda}\times }
 \times
 \left(\lambda_I-(I+1)\beta\right)
=\sum_i (\lambda_i^2-2i\lambda_i\beta).
\end{split}
\end{align}
\end{prp}

\begin{proof}
As in the proof of Proposition \ref{prp:young:1},
set 
$\lambda=(
 \stackrel{j_1}{\overbrace{n_1,\ldots,n_1}},
 \stackrel{j_2}{\overbrace{n_2,\ldots,n_2}},\ldots,
 \stackrel{j_l}{\overbrace{n_l,\ldots,n_l}})$
and $m_k\seteq j_1+\cdots+j_k$ ($k=1,\ldots,l$).
We can write the left hand side of \eqref{eq:prp:young:2} as
\begin{align*}
&F_2(\{m_k\},\{n_k\},\beta)
= \sum_{k=1}^l F_{2,k}(\{m_k\},\{n_k\},\beta),
\\
&F_{2,k}(\{m_k\},\{n_k\},\beta)
\seteq  (n_k-(m_k+1)\beta) (m_k-m_{k-1})(n_k-n_{k+1})
\\
&\times
 \prod_{i=1}^{k-1}
 \dfrac{(n_i-n_k)+\beta(m_k-m_{i-1})}{(n_i-n_k)+\beta(m_k-m_i)}
 \prod_{j=k+1}^{l}
 \dfrac{(n_k-n_{j+1})+\beta(m_j-m_k)}{(n_k-n_j)+\beta(m_j-m_k)}.
\end{align*}
The residues of $F_2$ are the same as those of $F_1$, and
by the similar calculation as in Proposition \ref{prp:young:1}, 
one can find that $F_2$ is a polynomial of $\beta$.
The behaviour of $F_2$ in the limit $\beta\to\infty$ shows that 
$F_2$ is a linear function of $\beta$.

Using the original expression \eqref{eq:prp:young:2},
we find that 
\[F_2(\{m_k\},\{n_k\},0)=\sum_{i}\lambda_i^2.\]
In order to determine the coefficient of $\beta$ in $F_2$, 
we rewrite $F_2$ as the rational function of $\beta^{-1}$, 
and take the limit $\beta^{-1}\to\infty$. 
The result is
\begin{align*}
&\lim_{\beta\to\infty} \big(\beta^{-1} F_2(\{m_k\},\{n_k\},\beta)\big)
\\
&=-\sum_{k=1}^l
   (m_k+1)(m_k-m_{k-1})(n_k-n_{k+1})\dfrac{m_k}{m_k-m_{k-1}}
\\
&=-\sum_{k=1}^l n_k(m_k-m_{k-1})(m_k+m_{k-1}+1).
\end{align*}
A moment thought shows that this becomes $-\sum_i 2 i \lambda_i$ 
if $\{m_k\}$ and $\{n_k\}$ correspond to $\lambda$.
Thus the proof is completed.
\end{proof}

\section{Conclusion and Remarks}\label{sec:rmk}

We have investigated the expansions of  Whittaker vectors 
for the Virasoro algebra in terms of Jack symmetric functions.
As we have mentioned in \S \ref{sect:intro}, 
the paper \cite[(3.18)]{AY:2009} proposed a conjecture on 
the factored expression for the Gaiotto state of the deformed Virasoro algebra.
using Macdonald symmetric functions.
However, our proof cannot be applied to this deformed case.
The main obstruction is that 
the zero-mode $T_0$ of the generating field $T(z)$ 
of the deformed Virasoro algebra behaves badly, 
so that one cannot analyse its action on Macdonald symmetric functions,
and cannot obtain a recursive formula similar to the one 
in Proposition \ref{prp:rec}.

It is also valuable to consider the $\calW(\frksl_n)$-algebra case.
In \cite{T:2009} a degenerate Whittaker vector is expressed in terms of 
the contravariant form of the $\calW(\frksl_3)$-algebra.
At this moment, however, we don't know how to treat 
Whittaker vectors for $\calW(\frksl_n)$-algebra.
It seems to be related to the higher rank analogues of 
the AGT conjecture (see \cite{W:2009} for examples).

\appendix
\section{AGT relation}\label{sec:AGT}

This appendix is devoted to the explanation of the AGT relation 
for pure $\SU(2)$ gauge theory, 
and its connection to the formula given in our main theorem.
This section is not necessary for the main argument of this paper.

\subsection{AGT relation for pure $\SU(2)$ gauge theory}

The original AGT conjecture \cite{AGT:2010} states the equivalence between 
the Liouville conformal blocks and the Nekrasov partition 
functions \cite{N:2003}.
In \cite{G:2009} the degenerated versions of the conjecture were proposed .
As the most simplified case, 
it was conjectured that the norm of the Gaiotto state of 
Virasoro algebra coincides 
with the Nekrasov partition function 
for the four-dimensional pure $\SU(2)$ gauge theory.

First we introduce the contravariant form (Shapovalov form) 
on the Verma module $M_h$.
Let us denote the (restricted) dual Verma module by $M^*_{h}$.
This is a right $\Vir_{c}$-representation generated by $\bra{h}$ with 
 $\bra{h}\Vir_{c,+}=0$ and $\bra{h} L_0=h\bra{h}$. 
The contravariant form is the bilinear map
\begin{align*}
\cdot: M^*_h \times M_h \to \bbC
\end{align*} 
determined by 
\begin{align*}
\bra{h}\cdot\ket{h}=1,\quad 
0\cdot \ket{h}=\bra{h}\cdot 0=0,\quad 
\bra{h}u\cdot\ket{h}=\bra{h}\cdot u\ket{h}\ (u\in \Vir).
\end{align*}

Fix a complex number $\Lambda$
\footnote{In this subsection we use the notations in the physics literatures.
Do not confuse this parameter $\Lambda$ and 
the notation $\Lambda$ of the ring of symmetric functions.}.
In this section we denote by $\ket{G}\in M_h$ 
\footnote{Do not confuse this symbol  $\ket{G}$ for the Gaiotto state and the symbol $\ket{h}$ for the highest weight vector.}
the Gaiotto state 
\begin{align*}
L_1\ket{G}=\Lambda^2 \ket{G},\quad 
L_n\ket{G}=0\ (n>1).
\end{align*}
normalised as 
\begin{align*}
\ket{G}=\ket{h}+\cdots.
\end{align*}
This normalisation condition means that 
the homogeneous component of $\ket{G}$ in $M_{h,0}$ is $\ket{h}$,
i.e.,
the coefficient $c_{\emptyset}(\alpha,\beta)$ in \eqref{eq:thm:expand} 
is set to be one.

Let us also define the anti-homomorphism 
\begin{align*}
\dagger: U(\Vir_{c,-}) \to U(\Vir_{c,+}), \quad L_{-n}\mapsto L_n.
\end{align*}
We will also denote the action of this map as $L_{-n}^\dagger=L_n$.
It induces a linear map 
$M_h\to M^*_h$, which is also written by $\dagger$. 
We define $\bra{G}\seteq (\ket{G})^\dagger$.

Next we recall the Nekrasov partition function 
(see \cite{N:2003} and \cite{B:2004,BE:2006}, \cite{NO:2006}, \cite{NY:2005}).
It has a geometric meaning, 
but here we only give  the next combinatorial expression. 
Let $r\in\bbZ_{\ge2}$ and 
$x,\ep_1,\ep_2,\vect{a}=(a_1,\ldots,a_r)$ be indeterminates.
Then the Nekrasov partition function 
$Z^{\rk=r}(x;\ep_1,\ep_2,\vect{a})$ 
for pure $\SU(r)$ gauge theory is defined to be:
\begin{align}
\label{eq:Nekrasov}
\begin{split}
Z^{\rk=r}(x;\ep_1,\ep_2,\vect{a})
\seteq
&\sum_{\vect{Y}}
\dfrac{ x^{|\vect{Y}|} }
{\prod_{1\le\alpha,\beta\le r} 
 n_{\alpha,\beta}^{\vect{Y}}(\ep_1,\ep_2,\vect{a})},
\\
n_{\alpha,\beta}^{\vect{Y}}(\ep_1,\ep_2,\vect{a})
\seteq
&\prod_{\square\in Y_\alpha}
 [-\ell_{Y_\beta}(\square)\ep_1+(a_{Y_\alpha}(\square)+1)\ep_2
 +a_\beta-a_\alpha]
\\
\times
&\prod_{\blacksquare\in Y_\beta}
 [(\ell_{Y_\alpha}(\blacksquare)+1)\ep_1- a_{Y_\beta}(\blacksquare)\ep_2
 +a_\beta-a_\alpha].
\end{split}
\end{align}
Here $\vect{Y}=(Y_1,\ldots,Y_r)$ is a $r$-tuple of partitions, 
$|\vect{Y}|\seteq |Y_1|+\ldots+|Y_r|$,
and $a_{Y}(\square)$, $\ell_{Y}(\square)$  are 
the arm and leg of the box $\square$ with respect to $Y$ as 
\begin{align*}
a_\lambda(\square)    \seteq \lambda_i-j,\quad 
\ell_\lambda(\square) \seteq \lambda'_j-i. 
\end{align*}
Note that for the case $i>\ell(\lambda)$ 
the number $\lambda_i$ should be taken as $\lambda_i=0$,
and  for $j>\lambda_1$ the number $\lambda'_j$ taken as  $\lambda'_j=0$. 
Thus $a_\lambda(\square)$ and $\ell_\lambda(\square)$ could be minus 
in general, 
although such cases don't occur in the norm of Jack symmetric functions.

Now the statement of the simplest case of the Gaiotto conjectures is 
\begin{align}\label{eq:AGT}
\left<G|G\right> \stackrel{?}{=} Z^{\rk=2}(x;\ep_1,\ep_2,\vect{a}).
\end{align}
Here the parameters are related as in Table \ref{table:VirNek}.
\begin{table}[htbp]
\centering
\begin{tabular}[c]{c c}
Virasoro & Nekrasov
\\
\hline
$c$ & $13+6(\ep_1/\ep_2+\ep_2/\ep_1)$
\\
$h$ & $\big((\ep_1+\ep_2)^2-(a_2-a_1)^2\big)/4 \ep_1 \ep_2$
\\
$\Lambda$ & $x^{1/4}/(\ep_1\ep_2)^{1/2}$
\end{tabular}
\\
\vskip 1em
\caption{Parameter correspondence}
\label{table:VirNek}
\end{table}

Note that this degenerate version of the AGT conjecture 
is proved by the method of Zamolodchikov-type recursive formula 
in the papers \cite{FL:2010} and \cite{HJS:2010}.

\subsection{Comparison of the inner products}

Our formula \eqref{eq:thm:expand} describes the Gaiotto state 
$\ket{G}$ by Jack symmetric functions.
In order to calculate the norm $\left<G|G\right>$,
we should compare the contravariant form $\cdot:M_h^*\to M_h$ and 
the inner product $\left<\cdot,\cdot\right>$ on $\Lambda$.

Let us recall the isomorphism $\iota_\beta$: 
\begin{align*}
\iota_\beta : \calF_\alpha\to\Lambda_\bbC,\quad
a_{-n}\ket{\alpha} \mapsto p_n\cdot \beta^{1/2}/\sqrt{2 }\ (n>0).
\end{align*}
In order to give the consistency between 
the bilinear form $\cdot:\calF_\alpha^*\times\calF_\alpha\to\bbC$ 
on the Heisenberg Fock space 
and the contravariant form $\cdot:M_h^*\times M_h\to\bbC$ 
on the Verma module of the Virasoro algebra,
we need to give the anti-homomorphism 
\begin{align*}
\omega:\calF_\alpha\to\calF_\alpha^*
\end{align*} 
so that 
\begin{align*}
\bra{h}u_1^\dagger \cdot u_2\ket{h}=
{}_\calF\bra{\alpha} \omega(\varphi(u_1)) \cdot \varphi(u_2) \ket{\alpha}_\calF
\end{align*}
holds for any $u_1,u_2\in U(\Vir_{c,+})$, 
where $\varphi: U(\Vir_c)\to \widehat{U}(\calH)$ 
is the bosonization map \eqref{eq:FF:phi}.

The consistent definition of $\omega$ is given as follows 
\cite{KM:1988}, \cite{TK:1986}:
\begin{align*}
\omega(a_n)=a_{-n}-2\rho\delta_{n,0},\quad
\omega(\rho)=-\rho.
\end{align*}
It implies for the parametrisation
$\rho=-(\beta^{1/2}-\beta^{-1/2})/\sqrt{2}$ that
$\omega(\beta^{1/2})=-\beta^{1/2}$.
Then we can spell out the inner product on $\Lambda$ 
which is consistent with the contravariant form on the Verma module $M_h$:
\begin{align*}
\left<p_n,p_m\right>=
{}_\calF\bra{\alpha} 
\omega(\sqrt{2}\beta^{-1/2}\, a_{-n})\cdot
\sqrt{2}\beta^{-1/2}\, a_{-m}
\ket{\alpha}_\calF
=-2n/\beta\cdot\delta_{n,m}.
\end{align*}
This is the inner product $\left<\cdot,\cdot\right>_{-2/\beta}$ 
defined in \eqref{eq:inner}.
But the Jack symmetric function orthogonal with respect to it 
is $P_\lambda^{(-2/\beta)}$, 
not $P_\lambda^{(1/\beta)}$ which is used in our expansion.

Thus the AGT relation \eqref{eq:AGT} is equivalent to 
\begin{align}\label{eq:AGT:2}
\sum_{\lambda,\mu\in\calP}\Lambda^{2|\mu|+2|\lambda|}
c_\lambda(\alpha,\beta)c_\mu(\alpha,\beta)
\prd{P_\lambda^{(1/\beta)},P_\mu^{(1/\beta)}}_{-2/\beta}
\stackrel{?}{=}Z^{\rk=2}(x;\ep_1,\ep_2,\vect{a}).
\end{align}
Now one may easily find that 
\begin{align*}
\prd{P_\lambda^{(1/\beta)},P_\mu^{(1/\beta)}}_{-2/\beta}
=0\quad\text{unless}\quad |\lambda|=|\mu|.
\end{align*}
Using this fact and 
comparing the homogeneous parts 
(the coefficients of $\Lambda^{4d}$ and those of $x^d$)
of both sides in \eqref{eq:AGT:2},
one finds that \eqref{eq:AGT} is equivalent to 
\begin{align}
\label{eq:AGT:3}
\begin{split}
\sum_{\substack{\lambda,\mu\,\vdash d}}
&c_\lambda(\alpha,\beta)c_\mu(\alpha,\beta)
\prd{P_\lambda^{(1/\beta)},P_\mu^{(1/\beta)}}_{-2/\beta}
\\
&\stackrel{?}{=}
(\ep_1 \ep_2)^{2d} 
\sum_{\substack{\lambda,\mu \in\calP \\ |\lambda|+|\mu|=d}}
\dfrac{1}{\prod_{1\le\alpha,\beta\le 2} 
 n_{\alpha,\beta}^{(\lambda,\mu)}(\ep_1,\ep_2,\vect{a})}
\end{split}
\end{align}
for each $d\in\bbZ_{\ge0}$.
In the right hand side we changed the notation $\vect{Y}\in\calP^2$ 
to the pair $(\lambda,\mu)\in\calP^2$.
Note that the ranges of running indexes in the left and right sides 
are different.
The equation \eqref{eq:AGT:3} seems to contain non-trivial relations 
among the `non-diagonal' pairings 
$\prd{P_\lambda^{(1/\beta)},P_\mu^{(1/\beta)}}_{-2/\beta}$.
According to the computer experiment,
these parings have complicated looks (in particular, no factored expressions) 
in general,
although each summand in the right hand side of \eqref{eq:AGT:3} is factored.
A combinatorial proof of \eqref{eq:AGT:3} would be another justification 
of the AGT relation \eqref{eq:AGT},
but we have no clue to show it directly at this moment.

We have another combinatorial restatement of \eqref{eq:AGT}.
If $\lambda\vdash d$, then
one can expand $P^{(1/\beta)}_\lambda\in\Lambda_\bbC^{d}$ 
by the basis $\{P^{(-2/\beta)}_\nu \mid \nu\vdash d\}$ of $\Lambda_\bbC^{d}$.
Let us express it as
\begin{align}\label{eq:apdx:exp}
P^{(1/\beta)}_\lambda
=\sum_{\nu\,\vdash d}\gamma^\nu_\lambda(\beta)P^{(-2/\beta)}_\nu,\quad
\gamma^\nu_\lambda(\beta)\in\bbC.
\end{align}
Then by an elementary calculation one finds 
that \eqref{eq:AGT} is equivalent to
\begin{align}
\label{eq:AGT:4}
\begin{split}
\sum_{\substack{\lambda,\mu,\nu\,\vdash d}}
&c_\lambda(\alpha,\beta)c_\mu(\alpha,\beta)
 \gamma^\nu_\lambda(\beta) \gamma^\nu_\mu(\beta) N_{\nu}(-2/\beta)
\\
&\stackrel{?}{=}(\ep_1 \ep_2)^{2d} 
\sum_{\substack{\lambda,\mu \in\calP \\ |\lambda|+|\mu|=d}}
\dfrac{1}{\prod_{1\le\alpha,\beta\le 2} 
 n_{\alpha,\beta}^{(\lambda,\mu)}(\ep_1,\ep_2,\vect{a})}.
\end{split}
\end{align}
Here we used the norm of Jack symmetric function
\begin{align*}
N_\nu(b)\seteq\prd{P_\nu^{(b)},P_\nu^{(b)}}_b=
 \prod_{\square\in\nu}
 \dfrac{a_\nu(\square)+b \ell_\nu(\square)+1}
       {a_\nu(\square)+b \ell_\nu(\square)+b}.
\end{align*}
According to the computer experiment,
the coefficient $\gamma_\lambda^\mu(\beta)$ 
in the expansion \eqref{eq:apdx:exp} 
doesn't have a factored expression in general,  
although it looks a little simpler 
than the pairing $\prd{P_\lambda^{(1/\beta)},P_\mu^{(1/\beta)}}_{-2/\beta}$.
One might 
find an explicit formula of $\gamma_\lambda^\mu(\beta)$.
However a direct proof of \eqref{eq:AGT:4} will require a manipulation on 
the changes of indexes from $\lambda,\mu,\nu\vdash d$ to 
$\lambda,\mu\in\calP$ with $|\lambda|+|\mu|=d$,
which seems to be hard at this moment.

\begin{ack}
The author is supported by JSPS Fellowships for Young Scientists (No.21-2241).
He expresses gratitude to the adviser Professor K\={o}ta Yoshioka 
and to Professor Yasuhiko Yamada for valuable discussion.
He would also like to thank the referees 
for their substantial suggestions on the improvements of the manuscript. 
\end{ack}


\end{document}